\documentclass{amsart}
\usepackage{amsmath,amssymb,dsfont,amsthm}
\usepackage{graphics}
\usepackage[all]{xy}
\usepackage{latexsym}
\usepackage{amsmath}
\usepackage{amsfonts}
\usepackage{amssymb}
\usepackage{tensor}
\usepackage{graphicx}

\date{\today}

\title[R torsion and analytic torsion for a conical frustum]{R torsion and analytic torsion of a conical frustum}
\thanks{2000 {\em Mathematics Subject Classification: 57Q10, 58J52}.}

\author{L. Hartmann  and M. Spreafico}

\address[Luiz Hartmann]{\tt UFSCar, Universidade Federal de S\~{a}o Carlos,  Brazil.  Partially supported by FAPESP 2010/16660-1}
\email{hartmann@dm.ufscar.br}

\address[Mauro Spreafico]{\tt ICMC, Universidade S\~{a}o Paulo, Brazil. Partially supported by 
FAPESP 2008/57607-6}
\email{mauros@icmc.usp.br}
\newtheorem{theo}{Theorem}
\newtheorem{lem}{Lemma}

\newtheorem{prop}{Proposition}

\renewcommand{\Im}{{\rm Im}}

\newcommand{\Sp}{{\rm Sp}}
\newcommand{\beq}{\begin{equation}}
\newcommand{\eeq}{\end{equation}}

\newcommand{\Z}{{\mathds{Z}}}
\newcommand{\R}{{\mathds{R}}}

\newcommand{\F}{{\mathds{F}}}
\newcommand{\Q}{{\mathds{Q}}}
\newcommand{\T}{{\mathcal{T}}}

\newcommand\e{{\rm e}}
\newcommand{\an}{{\rm{F}}}
\newcommand{\fr}{{\rm{F}}}
\newcommand{\A}{{\mathcal{A}}}
\renewcommand{\H}{{\mathcal{H}}}
\newcommand{\G}{{\Gamma}}
\renewcommand{\P}{{\mathcal{P}}}
\newcommand{\B}{{\mathcal{B}}}

\renewcommand{\b}{{\partial}}

\renewcommand{\det}{{\rm det}}

\newcommand{\CC}{{\mathcal{C}}}
\renewcommand{\v}{\mathrm{v}}
\newcommand{\semicolon}{,}%{\_}
\newcommand{\w}{\mathrm{w}}
\newcommand{\x}{\mathrm{x}}

\date{}

\DeclareMathOperator*{\Rz}{Res_0}

\DeclareMathOperator*{\Ru}{Res_1}

\begin{document}

\maketitle

\section{Introduction} 

Recently important advances have been made in the description of the analytic torsion of compact Riemannian manifolds with boundary \cite{Luc} \cite{DF} \cite{BM1} \cite{BM2}. In particular, in the last two works a formula for the analytic torsion of a compact oriented Riemannian manifold with boundary and absolute or relative boundary conditions was given. On the other side, in a series of works \cite{MS} \cite{HS1} \cite{HMS} \cite{HS2}, we presented explicit calculations of  the analytic torsion of some class of manifolds and pseudomanifolds. In particular formulas for the torsion of a cone  over a compact manifold where given. When working with cones, a natural question arises: if we truncate the cone we end up with a manifold. Does the analytic torsion of the cone coincide with some limit of the torsion of the truncated cone? This question was suggested to us by W. M\"{u}ller, and the answer is given in the Section \ref{lim} below, for an odd dimensional section. It turns out that some regularization is necessary, before taking the limit. The divergent terms are topological, in the sense that they come from the R torsion part of the analytic torsion, and non from the boundary term. More precisely, they come from the homology. Is then interesting to observe that the limit of the R torsion after this regularization process coincides with the intersection torsion, as expected, since the  cone is in general a pseudomanifold.

The natural way to tackle this question is to compute the Reidemeister torsion of the conical frustum and hence to apply  the Cheeger \cite{Che1} M\"{u}ller  \cite{Mul} theorem for manifolds with boundary in order to obtain the analytic torsion. Shortly, this means to  add the boundary contribution, as described in works of Br\"{u}ning and Ma \cite{BM1} \cite{BM2}. The last step will be to consider the limit case. This is the aim of this note, and is presented in the next four sections. In the last section we describe in details a particular case, namely the frustum over a circle. We present a detailed analysis of this case, that helps in understanding the general process. We also give an explicit calculation of the analytic torsion applying the definition. It is clear that the technique used for the circle admits a straightforward generalization to the case of any section; indeed, explicit calculations are given in the last section of \cite{HS2}, where however mixed boundary conditions were considered. Also note that the cone over a sphere is a  manifold, and hence in this case intersection torsion is replaced by genuine torsion: all the spaces involved in the limit process, namely the frustum and the cone, are regular manifolds; however a regularization is still necessary, since the homology is not trivial. 

\section{Geometric setting} Let $W$ be a compact connected oriented $m$ dimensional Riemannian manifold with metric $g$. The conical frustum (or truncated cone) over $W$ is the product manifold $\an W=[l_1,l_2]\times W$, where $0<l_1<l_2$, with metric (in the local coordinates $(x,y)$, where $y$ is a local system on $W$)
\[
g_\an=dx\otimes dx+x^2 g.
\]

The boundary of $\fr W$ is the disjoint union $\b \fr W=W_1\sqcup W_2$ of two copies of the $W$, $W_j=(W, l_j^2 g)$, with metric $l_j^2 g$. 

In order to deal with the R torsion of $\fr W$, we need some results on harmonic forms. By Hodge theory, it is clear that $\H^q(\fr W)$ is isomorphic to $\H^q(W)$. We need an explicit map. For we study the harmonic forms using the approach of \cite{Che0} (see also \cite{Nag}). It is clear that the formal solutions of the eigenvalues equation on $\fr W$ and on the cone over $W$ are the same, hence the next lemma follows (see \cite{HS2} Sections 3.3 and 8.1 for more details and for the notation).

\begin{lem}\label{l1} Let $\{\varphi_{{\rm har}}^{(q)},\varphi_{{\rm cl},n}^{(q)},\varphi_{{\rm ccl},n}^{(q)}\}$ be an orthonormal  base of $\Gamma(W,\Lambda^{(q)}T^* W)$
consisting of harmonic,  closed and coclosed  eigenforms of $\tilde\Delta^{(q)}$ on $W$. Let $\lambda_{q,n}$
denotes the  eigenvalue of $\varphi_{{\rm ccl},n}^{(q)}$ and $m_{{\rm ccl},q,n}$ its multiplicity. Define $\alpha_q =\frac{1}{2}(1+2q-m)$, $\mu_{q,n} = \sqrt{\lambda_{q,n}+\alpha_q^2}$, and $a_{\pm, q,n}=\alpha_q\pm \mu_{q,n}$. Then, 
all the solutions of the harmonic equation $\Delta u=0$, are convergent sums of forms of the following four types (as usual, a tilde denotes operations and quantities relative to the section):
\begin{align*}
\psi^{(q)}_{\pm, 1,n} =& x^{a_{\pm,q,n}} \varphi_{{\rm ccl},n}^{(q)},\\
\psi^{(q)}_{\pm, 2,n} =& x^{a_{\pm,q-1,n}} \tilde d\varphi_{{\rm ccl},n}^{(q-1)}+a_{\pm,q-1,n}x^{a_{\pm,q-1,n}-1}dx\wedge \varphi_{{\rm ccl},n}^{(q-1)},\\
\psi^{(q)}_{\pm,3,n} =& x^{a_{\pm,q-1,n}+2} \tilde d\varphi_{{\rm ccl},n}^{(q-1)}+a_{\mp,q-1,n}x^{a_{\pm,q-1,n}+1}dx\wedge \varphi_{{\rm ccl},n}^{(q-1)},\\
\psi^{(q)}_{\pm,4,n} =& x^{a_{\pm,q-2,n}+1}dx\wedge \tilde d\varphi_{{\rm ccl},n}^{(q-2)}.
\end{align*}

\end{lem}

Next, introducing absolute BC (as defined in \cite[Section 4]{RS} , or see \cite[Section 2]{HS1}): $\B_{\rm abs}(\omega)=0$ if and only if $\omega_{norm}|_{\b W}=0$ and $(d\omega)_{norm}|_{\b W}=0$, we have the following result, whose  proof is by direct verification: namely take the four types of forms as given in Lemma \ref{l1} and apply absolute BC to each. The unique forms that satisfy the absolute BC are the $\psi_{-,1,0}^{(q)}$, where $\lambda_{q,0}=0$ by definition. Sufficiency is easily verified: for $(\psi_{-,1,0}^{(q)})_{norm}=0$, and $(d \psi_{-,1,0}^{(q)})_{norm}=a_{\pm,q,n}x^{a_{-,q,n}-1}\psi_{-,1,0}^{(q)}$, and this vanishes at $x=l_1$ and $x=l_2$ if and only if $a_{-,q,n}=q-\sqrt{\lambda_{q,n}+q^2}=0$.  The result for relative BC is similar.

\begin{lem} \label{l2} The space of harmonic forms $\H^q_{\rm abs}(\fr W)$ coincides with the constant normal extension  %(in the normal direction) 
of the forms in $\H^q(W)$. The map $\omega\mapsto (-1)^q x^{m-2q} dx\wedge \omega$ defines an isomorphism of $\H^q(W)$ onto $\H^{q+1}_{\rm rel} (\fr W)$.  
\end{lem}

%The analytic torsion of the frustum is
%\[
%\log T_{\rm abs}((\fr W,g_\fr);\rho)=\log \tau_{\rm R}((\fr W,g_\fr);\rho)+A_{\rm BM, abs}(\b \fr W),
%\]
%where $\rho$ is an orthogonal representation of the fundamental group, and where the anomaly boundary term $A_{\rm BM, abs}(\b \fr W)$ is described in Section ? below.

\section{R torsion} In this section we calculate the R torsion of the frustum. For we first review  some necessary notation. % about R torsion.

%We will give an example of this calculation in Section \ref{lastsec}. We will follow here an alternative approach. 

\subsection{}\label{s3.1} We recall briefly the definition of the torsion of a finite chain complex of finite dimensional $\F$-vectors spaces (where $\F$ is a field of characteristic $0$)

\centerline{
\xymatrix{
\CC:& C_m\ar[r]^{\b_m}&C_{m-1}\ar[r]^{\b_{m-1}}&\dots\ar[r]^{\b_2}&C_1\ar[r]^{\b_1}&C_0.
}
}

Let $Z_q=\ker \b_q$, $B_q=\Im\b_{q+1}$, and $H_q=Z_q/B_q$. We assume that preferred bases $c_q=\{c_{q,j}\}$ and $h_q=\{h_{q,j}\}$ are  given for $C_q$ and $H_q$, respectively, for each $q$. Let  $b_q=\{b_{q,j}\}$ be a set of independent vectors in $C_q$ with $\b_q (b_q)\not=0$, and let $z_q=\{z_{q,j}\}$ be a set of independent vectors in $Z_q$ with $p(z_{q,j})=h_{q,j}$. Then, considering the sequence

\centerline{
\xymatrix{
0\ar[r]&B_q\ar[r]&Z_q\ar[r]_p& H_q\ar[r]&0,
}}

\noindent a basis for $Z_q$ is given by the basis $\b_{q+1}(b_{q+1})$ of $B_q$ and the set $z_q$. We denote this basis by $\b_{q+1}(b_{q+1})\semicolon z_q$ (see \cite{Mil} for details). By the same argument, the sequence

\centerline{
\xymatrix{
0\ar[r]&Z_q\ar[r]&C_q\ar[r]_{\b_q}& B_{q-1}\ar[r]&0,
}}

\noindent determine the basis $\b_{q+1}(b_{q+1})\semicolon z_q\semicolon b_q$ of $C_q$. Let $(\b_{q+1}(b_{q+1})\semicolon z_q\semicolon b_q/c_q)$ denote the matrix of the change of basis. Then, the torsion of $\CC$ is the class 
\[
\tau(\CC;\mathrm{v})=\prod_{q=0}^n [\det(\b_{q+1}(b_{q+1})\semicolon z_q\semicolon b_q/c_q)]^{(-1)^q},
\]
in $\F^\times/\{\pm 1\}$. It is easy to see that the torsion is independent of the graded bases $b=\{b_q\}$ and on the lifts $z=\{z_q\}$, but depends on the graded homology basis $h=\{h_q\}$. More precisely, $\tau(\CC;h)$ depends on the volume element $\v =\otimes_{q=0}^m h_q$ in $\otimes_{q=0}^m \Lambda^{r_q} H_q$, where $r_q={\rm rk}H_q$ (see for example \cite{Mul2}), and this explain the notation.

Now recall that the cylinder of the complex $\CC$ is the mapping cylinder of the identity $id:\CC\to \CC$, i.e. the complex 
$C_q(Cyl(\CC))=C_q\oplus C_{q-1}\oplus C_q$ with boundary
\[
\ddot\b=\left(\begin{matrix}\b&1&0\\0&-\b&0\\0&-1&\b\end{matrix}\right). 
\]
 
A preferred basis for  $C_q(Cyl(\CC))$ is $\ddot c_q=\{c_{q,j}\oplus 0\oplus 0,0\oplus c_{q-1,k}\oplus 0,0\oplus0\oplus c_{q,l}\}$.  By construction, $Cyl(\CC)$ has an homology graded preferred basis, and therefore its Whitehead torsion is well defined. We denote the preferred  basis of $H_q(Cyl(\CC))$ by $\ddot h_q$, and we let $\ddot z_q$ denotes  a lift of cycles of $\ddot h_q$.  Now we have the decomposition $B_q(Cyl(\CC))=\Im \ddot \b_{q+1}=(\Im \b_{q+1}+\ker \b_q)\oplus \Im \b_q\oplus (\Im \b_{q+1}+\ker \b_q)$, and hence %a basis for $B_q(Cyl(\CC))$ is given by the set $\{b_{q,j}\oplus0\oplus 0, 0\oplus0\oplus b_{q,k}, 0\oplus b_{q-1, l}\oplus 0, z_{q,n}\oplus 0\oplus -z_{q,m}\}$, and 
a set of independent elements in $C_q(Cyl(\CC))$ with non trivial image is %projecting onto these basis elements is
$\ddot b_q=\{b_{q,j}\oplus0\oplus 0, 0\oplus0\oplus b_{q,k}, 0\oplus b_{q-1, l}\oplus 0, 0\oplus z_{q-1,i}\oplus 0, 0\oplus \b_{q}(b_{q,h})\oplus 0\}$, that we denote by:  $b_{q}\oplus0\oplus 0, 0\oplus0\oplus b_{q}, 0\oplus b_{q-1}\oplus 0, 0\oplus z_{q-1}\oplus 0, 0\oplus \b_{q}(b_{q})\oplus 0$. A basis for $Z_q(Cyl(\CC))$ is then $\ddot b_q\semicolon \ddot z_q$, and  a   basis for $C_q(Cyl(\CC))$ is $\ddot \b_{q+1}(\ddot b_{q+1})\semicolon \ddot z_q\semicolon \ddot b_q$. Reordering and simplifying, we obtain the new basis
\begin{align*}
\ddot \b_{q+1}(\ddot b_{q+1})\semicolon \ddot z_q\semicolon \ddot b_q=&\b_{q+1}(b_{q+1})\oplus0\oplus 0\semicolon \ddot z_{q}\semicolon b_{q}\oplus 0\oplus 0\semicolon \\
 &0\oplus \b_{q}(b_{q})\oplus 0 \semicolon 0 \oplus z_{q-1}\oplus 0\semicolon 0\oplus  b_{q-1}\oplus0\semicolon\\
 &0\oplus0\oplus\semicolon  \b_{q+1}(b_{q+1})\semicolon 0\oplus 0\oplus z_q
 \semicolon 0\oplus 0\oplus b_q.
 \end{align*}

This gives (with some care at the higher dimensions)
\begin{align*}
[\det(&\ddot \b_{q+1}(\ddot b_{q+1})\semicolon \ddot z_q\semicolon \ddot b_q/\ddot c_q)]\\
%&=[(\ddot h_q/h_q)][\det( \b_{q+1}( b_{q+1})\semicolon  z_q\semicolon  b_q/ c_q)]^2[\det( \b_{q}( b_{q})\semicolon  z_{q-1}\semicolon  b_{q-1}/ c_{q-1})],
&=[(i_{*,q}^{-1}(\ddot h_q)/h_q)][\det( \b_{q+1}( b_{q+1})\semicolon  z_q\semicolon  b_q/ c_q)]^2[\det( \b_{q}( b_{q})\semicolon  z_{q-1}\semicolon  b_{q-1}/ c_{q-1})],
\end{align*}
where $i:\CC\to Cyl(\CC)$ denotes the inclusion, 
and for the torsion
\[
\tau(Cyl(\CC);\hat\v)= \frac{\tau(\CC;\v)}{\det({i_*})}.
\]

\subsection{}\label{s3.2} Let $(K,L)$ be a pair of connected finite cell complexes of dimension $m$, and $(\tilde K,\tilde L)$ its universal covering complex pair, and identify the fundamental group $\pi=\pi_1(K)$ with the group of the covering transformations of $\tilde K$. Note that covering transformations are cellular. Let $\CC((\tilde K,\tilde L);\Z)$ be the chain complex of $(\tilde K,\tilde L)$ with integer coefficients. The action of the group of covering transformations makes each chain group $C_q((\tilde K,\tilde L);\Z)$ into a module over the group ring $\Z\pi$, and each of these modules is $\Z \pi$-free  and finitely generated with preferred basis given by the natural choice of the $q$-cells of $K-L$. Since $K$ is finite it follows that $\CC((\tilde K,\tilde L);\Z)$ is free and finitely generated over $\Z \pi$. We obtain a complex of free finitely generated modules over $\Z\pi$ that we denote by $\CC((K, L);\Z\pi)$. Let $\rho:\pi\to O(\F, k)$ be an orthogonal representation of the fundamental group on a $\F$-vector space $V$ of dimension $k$, and consider the twisted complex $\CC((K,L);V_\rho)=V\otimes_{\Z\pi} \CC(( K, L);\Z\pi)$. Then, the torsion of $(K,L)$ with respect to the representation $\rho$ is the class
\[
\tau((K,L);\rho, \v)=\tau(\CC((K,L);V_\rho);\v),
\]
of $\F^\times/\{\pm 1\}$.

Next, let $W$ be an $m$ dimensional orientable compact connected Riemannian manifold with metric $g$ and possible boundary $\b W$. The torsion of $W$ can be defined taking any smooth triangulation or cellular decomposition of $W$. Moreover, the volume element $\v$ can also be fixed by using the metric structure. More precisely, given a graded orthonormal basis $a_q$ for the space of harmonic forms $\H^q(W)$, either with absolute or relative BC, and applying the De Rham map (see  for example \cite{RS})
\[
\A^{\rm abs}_{ q}=(-1)^q\P^{-1} \A^{m+1-q}_{{\rm rel}} \star:\H^q_{\rm abs}(W)\to H_q(W),
\]
we obtain a preferred homology graded basis $h=\A(a)$, that fix the volume element $\w=\A(\alpha)$, where $\alpha$ is the volume element determined by $a$. This gives the R torsion of $W$, and the relative R torsion of $(W,\b W)$:
\[
\tau_{\rm R}((W,g);\rho)=\tau(\CC(W;V_\rho);\A(\alpha)), \hspace{10pt} \tau_{\rm R}((W, \b W,g);\rho)=\tau(\CC((W,\b W);V_\rho);\A(\alpha)).
\]

\subsection{} It is clear that $\CC(\fr W)=\CC(W\times I)=Cyl(\CC(W))$, however in order to compute the R torsion we need to give to each complex the graded homology basis induced by the geometry.  For  we have at least two approaches: first apply the definition, and second use the exact sequence of the pair $(\fr W, W_2)$. We start with the first approach, and we will sckech the second one at the end of the section.

 Let denote by $\star_g$ the Hodge operator in the metric $g$. Let $g_2=l_2^2 g$. It is clear that $\star_{g_2}=l_2^{m-2q} \star_g$. Forms on $\fr W$ decompose as $\omega=\omega_1+dx\wedge \omega_2$. Writing $\omega(x,y)=f_1(x)\omega_1(y)+f_2(x)dx\wedge \omega_2(y)$, a simple calculation gives 
%(see \cite{HS1} Section 2.1 and equation (3.1)) gives
\[
\star_{g_F} \omega=x^{m-2q+2}f_2(x)\star_g \omega_2(y)+(-1)^q x^{m-2q} f_1(x)\wedge \star_g \omega_1(y).
\]

As vector spaces $\H^q(W)=\H^q(W_2)$. If $\omega$ is in $\H^q( W)$, denote the the constant extension of $\omega$ in $\H^q_{\rm abs}(\fr W)$ by $\ddot\omega$. Then:
\[
|| \omega||_{l_2^2  g}^2 = \int_W a_{q,j} \wedge \star_{l_2^2  g} \omega  = l_2^{m-2q} \int_W  \omega \wedge \star_{g} \omega = l_2^{m-2q} ||\omega||^2_{ g}, 
\] 
and
\begin{align*}
|| \ddot\omega ||_{g_F}^2 &= \int_{FW} \ddot\omega \wedge \star_{g_F} \ddot\omega = \int_{l_1}^{l_2} x^{m-2q} dx \int_W  \omega \wedge \star_{ g} \omega= \Gamma_q ||\omega||_g^2,
\end{align*}
where
\[
\Gamma_q=\left\{
\begin{array}[c]{c}
\frac{1}{m+1-2q}\left(l_2^{m+1-2q}-l_1^{m+1-2q}\right),\qquad {\rm if}\qquad |m+1-2q|\neq0, \\
\ln\frac{l_2}{l_1}, \qquad\qquad\qquad\qquad\qquad\;\;\;\qquad {\rm if} \;\qquad m+1-2q=0.
\end{array}
\right.\\  
\]

By the definition of the De Rham map $\A^{\rm abs}_{g_F, q}=(-1)^q\P^{-1} \A^{m+1-q}_{g_F,{\rm rel}} \star_{g_F}$ on $\fr W$ (see \cite{RS} Section 3), we have the following commutative diagram of isometries of vectors spaces (where the $\dagger$ denotes the  dual block complex)

\centerline{
\xymatrix{\H^q_{\rm abs}(\fr W)\ar[r]^{\star_{g_F}}&\H^{m-q+1}_{\rm rel}(\fr W)\ar[r]^{\A_{g_F, {\rm rel}}^{m-q+1}}&H^{m-q+1}((\fr W)^\dagger, (\b \fr W)^\dagger)&H_q(\fr W)\ar[l]_{\hspace{40pt}\P}\\
\H^q( W)\ar[r]^{\star_{g}}\ar[d]^{\frac{1}{\sqrt{l_2^{m-2q}}}}\ar[u]_{\frac{\ddot{(\_)}}{\sqrt{\Gamma_q}}}&
\H^{m-q}( W)\ar[r]^{\A^{m-q}_{g}}\ar[u]_{\frac{(-1)^q x^{m-2q}dx\wedge \_}{\sqrt{\Gamma_q}}}\ar[d]^{\sqrt{l_2^{m-2q}}}&H^{m-q}( W^\dagger)\ar[u]\ar[d]&H_q( W)\ar[l]_{\P}\ar[u]_{(-1)^q\sqrt{\Gamma_q}}\ar[d]^{\sqrt{l_2^{m-2q}}}\\
\H^q( W_2)\ar[r]^{\star_{g_2}}&\H^{m-q}( W_2)\ar[r]^{\A_{g_2}^{m-q}}&H^{m-q}( W_2^\dagger)&H_q( W_2)\ar[l]_{\P}\\
}}

Commutativity of the first square follows by the given formula for the Hodge operator and Lemma \ref{l2}. Commutativity of the other squares follows by construction. For suppose a cell decomposition of $W$ is fixed. Then, a cell decomposition of $\fr W$ is determined with $q$-cells either the  $q$-cells of $W$ or the product $I\times c$, where $c$ is a $(q-1)$-cell of $W$. Namely, we are using the  direct sum decomposition of the cellular chain complex $C_q(\fr W)=C_q(W)\oplus C_{q-1}(W)\oplus C_q(W)$. This fix a preferred basis for the chain vector spaces. It is clear that the dual block in $\fr W$  of $\{0\}\times c$ is  $[0,1/2]\times c^\dagger$, where $c^\dagger $ is the dual block of $c$ in $W$. Next, recall  the De Rham maps  $\A_g^{q}$ and $\A^q_{g_F,{\rm rel}}$  on $W$ and on $\fr W$ are defined respectively by
\[
\A_g^q(\omega)(c^\dagger)=\int_{c^\dagger} \omega,\hspace{30pt} \A^q_{g_F,{\rm rel}}(\omega)(d^\dagger)=\int_{d^\dagger} \omega,
\]
where $c^\dagger$ denotes the dual block of a cell $c$ of $W$, and $d^\dagger$ the dual block of a cell $d$ in $\fr W-\b\fr W$.  It follows that $\A^{m-q+1}_{g_F,{\rm rel}}(\omega)(c^\dagger)$ is non vanishing only on the dual block $d^\dagger$ of $(q+1)$-cells $d$ of  type $\{0\}\times c$. Hence, $\omega$ must have a non trivial normal component, i.e. $\omega=f_2(x) dx\wedge \omega_2(y)$, and the unique contributions are
\[
 \A^q_{g_F,{\rm rel}}(\omega)(d^\dagger)=\int_{(c,0)^\dagger} f_2(x)dx\wedge \omega_2(y)=\int_{l_1}^{l_2} f_2(x) dx \int_{c^\dagger} \omega_2(y).
 \]

This gives the isomorphisms in the  vertical lines of the last square, and their coefficients. 

Now realize the imbedding of $W$ in $\fr W$ as $W_2$.  Let $\alpha_q$ be an orthonormal base for $\H^q(W)$. Then,  an orthonormal base for $\H^{q}(W_2)$ is $l_2^{-\frac{m-2q}{2} }\alpha_q$, and applying the de Rham maps   we obtain
\begin{align*}
\A_{g_2,q}(l_2^{-\frac{m-2q}{2}}\alpha_{q,j}) &= l_2^{-\frac{m-2q}{2}} \A_{g_2,q}(\alpha_{q,j}) = l_2^{-\frac{m-2q}{2}} \P_q^{-1} \A_{g_2}^{m-q} \star_{  g_2} (\alpha_{q,j})\\
&= l_2^{-\frac{m-2q}{2}} l_2^{m-2q}  \P_q^{-1} \A_g^{m-q} \star_g (\alpha_{q,j}) = l_2^{\frac{m-2q}{2}} \A_{ g,q}(\alpha_{q,j}),
\end{align*} 
and hence  a basis $h_q$ for $H_q(W_2)$ is $l_2^{\frac{m-2q}{2}} \A_{g,q}(\alpha_q)$, and $z_q=l_2^{\frac{m-2q}{2}} \A_{g,q}(\alpha_q)$. Next, consider the conical frustum $(FW,g_F)$. An orthonormal basis for $\H^{q}_{\rm abs}(\fr W)$ is $\G_q^{-\frac{1}{2}} \alpha_q$, and
\begin{align*}
\A^{\rm abs}_{ g_F,q}\left(\G_q^{-\frac{1}{2}}\alpha_{q,j}\right) &= \G_q^{-\frac{1}{2}} \A^{\rm abs}_{g_F,q}(\alpha_{q,j})
= \G_q^{-\frac{1}{2}} \P_q^{-1} \A^{m+1-q}_{g_F,{\rm rel}} \star_{g_F} (\alpha_{q,j})\\
&= \G_q^{-\frac{1}{2}} \G_q \P_q^{-1} \A^{m-q}_{g} \star_{ g} (\alpha_{q,j}) =\G_q^{\frac{1}{2}}\A_{ g,q}(\alpha_{q,j}).
\end{align*} 

Then the basis $\ddot h_q$ is  $\frac{\G_q^{\frac{1}{2}}}{l_2^{\frac{m-2q}{2}}} (h_q\oplus 0\oplus 0)$, and hence $\ddot z_q=\frac{\G_q^{\frac{1}{2}}}{l_2^{\frac{m-2q}{2}}} (z_q\oplus 0\oplus 0)$. This gives $\det (i_*)=\frac{l_2^{\frac{m-2q}{2}}}{\G_q^{\frac{1}{2}}}$, thus 
\begin{align*}
 [\det(\ddot \b_{q+1}(\ddot b_{q+1})\semicolon \ddot z_q\semicolon \ddot b_q/\ddot c_q)]=%&\left[\det( \b_{q+1}( b_{q+1})\semicolon \A^{\rm abs}_{ g_F,q}\left(\G_q^{-\frac{1}{2}}\alpha_{q}\right) \semicolon  b_q/ c_q)\right]\\
% & \times\left[\det( \b_{q+1}( b_{q+1})\semicolon  \A_{g_2,q}(l_2^{-\frac{m-2q}{2}}\alpha_{q})\semicolon  b_q/ c_q)\right]\\
%&\times \left[\det( \b_{q}( b_{q})\semicolon  \A_{g_2,q}(l_2^{-\frac{m-2q}{2}}\alpha_{q})\semicolon  b_{q-1}/ c_{q-1})\right]\\
& \left(\frac{\G_q}{l_2^{m-2q}}\right)^\frac{r_q}{2}[\det\rho( \b_{q+1}( b_{q+1})\semicolon  \A(\alpha_q)\semicolon  b_q/ c_q)]^2\\
&\times [\det\rho( \b_{q}( b_{q})\semicolon  \A(\alpha_{q-1})\semicolon  b_{q-1}/ c_{q-1})],
\end{align*}
and
\[
\tau_{\rm R}((\fr W,g_F);\rho)=\prod_{q=0}^{m} \left(\frac{\G_q}{l_2^{m-2q}}\right)^\frac{(-1)^qr_q}{2}[\det\rho( \b_{q+1}( b_{q+1})\semicolon  \A(\alpha_q)\semicolon  b_q/ c_q)]^{(-1)^q}.
\]

We have proved the following proposition.

\begin{prop}\label{p1} The R torsion of the conical frustum is:
\[
\log  \tau_{\rm R}((\fr W,g_F);\rho)=\log \tau_{\rm R}((W,l_2^2 g);\rho)+\log \tau(\T),
\]
where
\begin{align*}
\log\tau(\T)=&\frac{1}{2}\sum_{q=0}^{2p} (-1)^q r_q \log\frac{l_2^{2p+1-2q}-l_1^{2p+1-2q}}{(2p+1-2q)l_2^{2p-2q}},&m&=2p, p\geq 0,\\
\log\tau(\T)=&\frac{1}{2}\sum_{q=0,q\not= p}^{2p-1} (-1)^q r_q \log\frac{l_2^{2p-2q}-l_1^{2p-2q}}{(2p-2q)l_2^{2p-1-2q}}&m&=2p-1, p\geq 1,\\
&+\frac{(-1)^p}{2}r_p \log\log\frac{l_2}{l_1},
\end{align*}
where $r_q={\rm rk} H_q(W)$.
\end{prop}

We conclude this section with a second proof of Proposition \ref{p1}. Consider the short exact sequence of chain complexes associated to the pair $(\fr W, W_2)$,
\[
0\to \CC(W_2)\to \CC (\fr W)\to \CC(\fr W ,W_2)\to 0,
\]
by Milnor \cite[Section 3]{Mil}, we have  
\[
\log  \tau_{\rm R}((\fr W,g_F);\rho)=\log \tau_{\rm R}((W,l_2^2 g);\rho)+\log  \tau_{\rm R}((\fr W, W_2),g_F);\rho)+\log \tau(\T),
\]
where the complex $\T$ is defined by the long exact homology sequence of the pair, namely
\beq\label{sec}
\T:{\xymatrix@C=0.5cm{& \ldots \ar[r] &H_q(W_2)\ar[r] & \ar[r] H_q(\fr W) & \ar[r] H_q(\fr W,W_2)& \ldots}},
\eeq
with $\T_{3q+2} =  H_q (W_2)$, $\T_{3q+1} =  H_q (\fr W)$ and $\T_{3q} = H_q (\fr W, W_2)$. It is clear that both the relative torsion and the relative homology are trivial. Therefore the torsion of $\T$ is given by the graded product of the torsions of the isomorphisms: $i_{*,q}:H_q(W)\to  H_q(\fr W)$. Using the graded homology basis given above, we can now compute the determinants of the change of basis in the vector spaces of the sequence in equation (\ref{sec}). At
 $ H_q(W,W_2)$ the determinant is $1$, at $ H_q(W_2)$ is $1$ and at $ H_q(FW)$ is $\left(\frac{l_2^{m-2q}}{\G}\right)^{\frac{r_q}{2}}$, where $r_q$ is the rank of the homology. Applying the definition of Reidemeister torsion to the complex $\T$, we obtain (where $D$ denotes the determinant of the matrix of the change of basis)
\begin{align*}
\log\tau(\T) =& \sum_{q=0}^{3m} (-1)^q \log D(\T_{q}) \\
=&\sum_{q=0}^{m} (-1)^{3q} \log D(H_{q}(FW,W_2)) + \sum_{q=0}^{m} (-1)^{3q+1}\log D(H_q(FW))\\
&+ \sum_{q=0}^{m} (-1)^{3q+2}\log D(H_q(W_2))\\
=&\sum_{q=0}^{m} (-1)^{3q+1}\log D(H_q(FW))\\
=&\sum_{q=0}^{m} (-1)^{q+1}\frac{r_q}{2} \log \frac{l_2^{m-2q}}{\G_q}.
\end{align*}

\section{Analytic torsion} 
\label{s4}

Using the works of Br\"{u}ning and Ma \cite{BM1}  \cite{BM2}, the Cheeger M\"{u}ller theorem for an oriented compact connected Riemannian $n$-manifold $(M,g)$ with boundary reads (see \cite[Section 6]{HS1}  or \cite[Section 2.3]{HS2}  for details on our notation) 
\begin{align*}
\log T_{\rm abs}((M,g);\rho)&=\log\tau_{\rm R}((M,g);\rho)+\frac{{\rm rk}(\rho)}{4}\chi(\b M)\log 2+{\rm rk}(\rho)A_{\rm BM,abs}(\b M),\\
\log T_{\rm rel}((M,g);\rho)&=\log\tau_{\rm R}((M,\b M,g);\rho)+\frac{{\rm rk}(\rho)}{4}\chi(\b M)\log 2+{\rm rk}(\rho)A_{\rm BM,rel}(\b M),
\end{align*}
where $\rho$ is an orthogonal representation of the fundamental group, and where the boundary anomaly term of Br\"{u}ning and Ma is defined as follows. Using the notation of \cite{BM1} (see \cite[Section 2.2]{HS2} for more details) for $\Z/2$ graded algebras, we  identify an antisymmetric  endomorphism $\phi$ of a finite dimensional vector space $V$ (over a field of characteristic zero) with the element 
$\hat \phi=\frac{1}{2}\sum_{j,k=1}^n \langle\phi(v_j),v_k\rangle \hat v_j\wedge \hat v_k$,
of $\widehat{\Lambda^2 V}$. For the elements $\langle\phi(v_j),v_k\rangle$ are the entries of the tensor representing $\phi$ in the base $\{v_k\}$, and this is an antisymmetric matrix. Now assume that $r$ is an antisymmetric endomorphism of $\Lambda^2 V$. Then, $(R_{jk}=\langle r(v_j),v_k\rangle)$ is a tensor of two forms in $\Lambda^2 V$. We extend the above construction identifying $R$ with the element
\[
\hat R=\frac{1}{2}\sum_{j,k=1}^n \langle r(v_j),v_k\rangle\wedge \hat v_j\wedge \hat v_k,
\]
of $\Lambda^2 V\wedge \widehat{\Lambda^2 V}$. This can be generalized to higher dimensions. In particular,
all the construction can be done taking the dual $V^*$ instead of $V$. Accordingly to \cite{BM1}, we  define the following forms (where $i:\b M\to M$ denotes the inclusion)
\begin{align*}
\mathcal{S}&=\frac{1}{2}\sum_{k=1}^{n-1}(i^*\omega-i^*\omega_0)_{0 k}\wedge\hat e^*_{k}\\
\widehat {i^*\Omega}&=%\mathcal{R}^{T C_{\alpha}S^{m-1}_{l\sin\alpha}}|S^{m-1}_{la}=
\frac{1}{2}\sum_{k,l=1}^{n-1}i^*\Omega_{k l}\wedge\hat e^*_{k}\wedge  \hat e^*_{l},&
\hat{\Theta}&=\frac{1}{2}\sum_{k,l=1}^{n-1}\Theta_{kl} \wedge \hat  e^*_{k}\wedge  \hat e^*_{l}.
\end{align*}

Here, $\omega$ and $\omega_0$  are the connection one forms associated to the metrics $g_0$ and $g_1=g$, respectively,  where $g_0$ is a suitable deformation of $g$ that is a product near the boundary. $\Omega$ is  the curvature two form  of $g$, $ \Theta$ is the curvature two form of the boundary (with the metric induced by the inclusion), and $\{e_k\}_{k=0}^{n-1}$ is an orthonormal base of $T M$ (with respect to the metric $g$). Then, setting 
\[ 
B=\frac{1}{2}\int_0^1\int^B
\e^{-\frac{1}{2}\hat{\Theta}-u^2 \mathcal{S}^2}\sum_{k=1}^\infty \frac{1}{\Gamma\left(\frac{k}{2}+1\right)}u^{k-1}
\mathcal{S}^k du, 
\]
the anomaly boundary term is
\[
A_{\rm BM,abs}(\b M)=(-1)^{n+1}A_{\rm BM,rel}(\b M)=\frac{1}{2}\int_{\b M} B.
\]

It is not too difficult to see that in the case of the frustum $\fr W$, the boundary term $A_{\rm BM}(W_j)$  is independent on $l_j$, either with absolute or relative BC. For let  $\{b_k\}_{k=1}^m$ and $\{e_k\}_{k=0}^{m}$ denote   local orthonormal bases of $TW$ and $T \fr W$ respectively. Then, direct calculations (see \cite[Section 3.2]{HS2}, see also Section \ref{lastsec} for an example) give
\begin{align*}
\mathcal{S}&=-\frac{1}{2l}\sum_{k=1}^{m}e^*_k\wedge e^*_{k}=-\frac{l}{2}\sum_{k=1}^{m}b^*_k\wedge b^*_{k}=-\frac{1}{2}\sum_{k=1}^{m}b^*_k\wedge e^*_{k},
\end{align*}
and $\Theta_{jk}=\tilde\Omega_{jk}$. On the other side, it is also clear from the definition that the boundary terms on the two boundaries will have either opposite sign or the same sign depending on the dimension $m$ of the boundary. We have the following result.

\begin{lem} The anomaly boundary term on the frustum is: if $m$ is odd
\begin{align*}
A_{\rm BM,abs}(\b \fr W)=-A_{\rm BM,rel}(\b \fr W)&=0, &
 A_{{\rm BM,abs} W_1,{\rm  rel}  W_2 }(\b \fr W)&=\int_{W} B; 
\end{align*}
if $m$ is even
\begin{align*}
A_{\rm BM,abs}(\b \fr W)=A_{\rm BM,rel}(\b \fr W)=
 A_{{\rm BM,abs} W_1,{\rm  rel}  W_2 }(\b \fr W)&=\int_{W} B. 
\end{align*}
\end{lem}

%This proves that the analytic torsion with either absolute or relative BC is independent on the boundary, while that with mixed BC is pure boundary:

\begin{prop} The analytic torsion of the conical frustum is
\begin{align*}
\log T_{\rm abs}((\fr W,g_\fr);\rho)=&\log \tau_{\rm R}((\fr W,g_\fr);\rho)+\frac{1}{2}{\rm rk}(\rho)\chi(W)\log 2\\
&+\frac{1-(-1)^{m+1}}{2}{\rm rk}(\rho)\int_W B,\\
% \log T_{\rm rel}((\fr W,g_\fr);\rho)=&\log \tau_{\rm R}((\fr W,\b \fr W,g_\fr);\rho)+\frac{1}{4}{\rm rk}(\rho)\chi(W)\log 2\\
 %&+\frac{1-(-1)^{m+1}}{2}{\rm rk}(\rho)\int_W B,
\log T_{{\rm abs} W_1,{\rm  rel}  W_2 }((\fr W,g_\fr);\rho)=&\frac{1}{2}{\rm rk}(\rho)\chi(W)\log 2+{\rm rk}(\rho)\int_W B.
\end{align*}
%and by duality 
%\begin{align*}
%\log T_{\rm abs}((\fr W,g_\fr);\rho)&=(-1)^m \log T_{\rm rel}((\fr W,g_\fr);\rho),\\ 
%\log T_{{\rm abs} W_1,{\rm  rel}  W_2 }((\fr W,g_\fr);\rho)&=(-1)^m \log T_{{\rm abs} W_2,{\rm  rel}  W_1 }((\fr W,g_\fr);\rho).
%\end{align*}

\end{prop}

\section{Limit case}
\label{lim}

In this section we study the limit case $l_1\to 0^+$, and the relation with the torsion of the cone $C_{l_2} W=\fr W/W_1$. For we first give the formula for the analytic torsion of the cone (formulas for relative BC follow by duality, as proved in Theorems 1.1 and 1.2 of \cite{HS2}). In this section we analyze the case of odd dimensional section, so we assume $m=2p-1$, $p>1$; we also assume ${\rm rk}(\rho)=1$. %The formula for odd section was given and proved in \cite{HS1}. The formula for even section was nor given neither proved, however as observed in that paper, the proof is exactly the same as in the odd case, only the result is quite different.

\begin{theo}[\cite{HS2}]\label{t01} The analytic torsion on the cone $C_l W$ on an orientable compact connected Riemannian manifold $(W,g)$ of odd dimension $2p-1$ is
\begin{align*}
\log T_{\rm abs}((C_lW,g_C);\rho_0))= &\frac{1}{2} \sum_{q=0}^{p-1} (-1)^{q+1} {\rm rk}H_q(W;\Q)\log \frac{2(p-q)}{l}\\
&+\frac{1}{2} \log T((W,l^2g);\rho_0)+A_{\rm BM,abs}(\b C_l W).
\end{align*}

%The torsion when the dimension of $W$ is  $2p$ even is
%\begin{align*}
%\log T_{\rm abs}(C_lW)= &\sum_{q=0}^{p} (-1)^{q+1} r_q \left(q-p-\frac{1}{2}\right)\log l +\frac{1}{2}\sum^{p-1}_{q=0}(-1)^{q+1} A_{q,0,0}(0)\\
%&+ \frac{1}{2} \chi(W) \log2+\frac{1}{2}\sum_{q=0}^{p-1} (-1)^{q+1}r_q \left(\log(2p-2q+1) +2\log(2p-2q-1)!!\right)\\
%&+A_{\rm BM,abs}(\b C_l W).
%\end{align*}
%where ($\alpha_q$, $m_{{\rm ccl}, q,n}$ and  and $\mu_{q,n}$ were defined in Lemma \ref{l1})
%\[
%\A_{0,0,q}(s)= \sum_{n=1}^{\infty}
%\left(\log\left(1-\frac{\alpha_q}{\mu_{q,n}}\right) - \log\left(1+\frac{\alpha_q}{\mu_{q,n}}\right)\right)\frac{m_{{\rm ccl}, q,n}}%{\mu_{q,n}^{2s}}.
%\]

\end{theo}

A simple calculation (using for example the variational formula for the torsion) shows that
\[
\log\tau_{\rm R}((W,l^2 g);\rho)=\log\tau_{\rm R}((W, g);\rho)+\frac{1}{2}\sum_{q=0}^{m}(-1)^q r_q (m-2q)\log l,
\]
and, by duality
\[
\sum_{q=p}^{2p-1} (-1)^q r_q (2p-1-2q)=
\sum_{q=0}^{p-1}(-1)^q r_q (2p-1-2q)=\frac{1}{2}\sum_{q=0}^{m}(-1)^q r_q (m-2q).
\]

Thus the formula for the analytic torsion of the cone $C_{l_2} W$ reads:
\begin{align*}
\log T_{\rm abs}((C_{l_2}W,g_C);\rho_0)= &\frac{1}{2} \log \tau_{\rm R}((W,g);\rho_0)+\frac{1}{2} \sum_{q=0}^{p-1} (-1)^{q} r_q(2p-1-2q)\log l_2\\
&+\frac{1}{2} \sum_{q=0}^{p-1} (-1)^{q+1} r_q\log \frac{2(p-q)}{l_2}+A_{\rm BM, abs}(\b C_{l_2} W).
\end{align*}

Consider the formula for the R torsion of the frustum given in Proposition \ref{p1}. It is clear that in the limit $l_1\to 0^+$ the last $p$ terms diverge. This suggests the following approach. 

Let $(W,g)$ be an oriented compact connected Riemannian manifold of dimension $m$. Such a space has a class of distinguished CW decompositions (given by the  smooth triangulations). Let $K$ one of these CW decompositions, and let  $c$ the  preferred graded basis of the chain complex $\CC(W;V_\rho)$ given by the cells, as described in Section \ref{s3.2}. Fix the sets $b$ and $z$ as in Section \ref{s3.1}, and let $\A(\alpha)$ denote the graded basis for homology induced by the metric structure as in Section \ref{s3.2}, and use the notation
\[
D_q=[\det\rho( \b_{q+1}( b_{q+1})\semicolon  z_q\semicolon  b_q/ c_q)]\in \F^\times/\{\pm 1\}.
\]

Let $K^\dagger$ denote the dual block complex, and $K_{(q)}$ the $q$-skeleton. It is clear that the $p-1$ homology of $K_{(p-1)}$ coincides with the cycles of $K$, and the bijection is cellular. Consider the torsion
\[
\tau (W_{(p-1)};\rho,\w)=\tau (\CC (K_{(p-1)};V_\rho),\w)=\prod_{q=0}^{p-1} [\det\rho( \b_{q+1}( b_{q+1})\semicolon  z_q\semicolon  b_q/ c_q)]\in \F^\times/\{\pm 1\},
\]
where the $z_q$ are cycles projecting onto a  basis for $H_q(K_{(p-1)})$, and  
$\w$ is the induced volume element, as in Section \ref{s3.2}. We can fix $\w$ using the geometry. Since $K^\dagger$ is another decomposition of $W$, there is a common subdivision $T$ of $K$ and $K^\dagger$. The identity maps $id:W=|K|\to W=|T|$ and $id:W=|K^\dagger|\to W=|T|$ are cellular, and hence restrict to maps $id_{(q)}:K_{(q)}\to T_{(q)}$ and $id_{(q)}:(K^\dagger)_{(q)}\to T_{(q)}$, i.e. $T_{(q)}$ is a common subdivision of $K_{(q)}$ and $(K^\dagger)_{(q)}$. It follows that $K_{(q)}$ and $(K^\dagger)_{(q)}$ have the same torsion up to the choice of the homology volume elements, by \cite[7.1]{Mil}.  Consider the chain complex associated to $(K^\dagger)_{(p-1)}$. Let 
\[
D^\dagger_q=[\det\rho( \b^\dagger_{q+1}( b^\dagger_{q+1})\semicolon  z^\dagger_q\semicolon  b^\dagger _q/ c^\dagger_q)]\in \F^\times/\{\pm 1\}.
\]

By duality $h^\dagger_q=h_{2p-1-q}$, $D^\dagger_q=D^{-1}_{2p-1-q}$, and hence
\[
\prod_{q=0}^{p-1}\left( D^\dagger_q\right)^{(-1)^q}=\prod_{q=p}^{2p-1} D_q^{(-1)^q},
\]
and
\[
\prod_{q=0}^{2p-1} D_q^{(-1)^q}=\prod_{q=0}^{p-1}D_q^{(-1)^q}\prod_{q=0}^{p-1}\left( D^\dagger_q\right)^{(-1)^q}.
\]

It is clear that the basis $\A(\alpha_q)$ gives an homology basis for $H_q(K_{(p-1)})$ for all $q<p-1$. Moreover, $z_{p-1}=\b_p(b_{p})\semicolon \A(\alpha_{p-1})$ gives a basis for $H_{p-1}(K_{(p-1)})$. This basis depends on the   $b_{p}$, however, if we change the set $b_{p}$ by $b'_{p}$, we have
\[
D'_{p-1}=k D_{p-1},
\]
for some  field unit $k$ (up to sign). Also, the dual basis change gives the change
\[
(D^\dagger)'_{p-1}=k^{-1}D^\dagger _{p-1}.
\]

It follows that there exists a family $\B$ of homology basis  of $H_{p-1}(K_{(p-1)})$, but a unique volume element $\x$, such that 
\[
1=\frac{\tau (W_{(p-1)};\rho,\x)}{\tau ((W^\dagger)_{(p-1)};\rho,\x^\dagger)}.
\]

This fix the volume element $\w=\x$, and with this choice
\[
\tau_{\rm R}((W,g);\rho)=(\tau (W_{(p-1)};\rho,\x))^2.
\]
%We can therefore fix a  unique field unit 
%\[
%x=\sqrt\frac{\tau (K_{(p-1)};\rho,\w)}{\tau ((K^\dagger)_{(p-1)};\rho,\w^\dagger)},
%\]
%and fix the volume element $\w$ by requiring that 

Back to the frustum, we have that
\begin{align*}
\tau_{\rm R}((FW,g_F);\rho)&=\prod_{q=0}^{m}  \left(\frac{\G_q}{l_2^{m-2q}}\right)^{(-1)^q\frac{ r_q}{2}}D_q^{(-1)^q} \\
&= \prod_{q=0}^{m}  \left(\frac{\G_q}{l_2^{m-2q}}\right)^{(-1)^q\frac{ r_q}{2}}
[\det\rho( \b_{q+1}( b_{q+1})\semicolon  \zeta_q\semicolon  b_q/ c_q)]\\
&= \left(\prod_{q=0}^{m}  \left(\frac{\G_q}{l_2^{m-2q}}\right)^{(-1)^q\frac{ r_q}{2}}\right)
\left(\tau (W_{(p-1)};\rho,\x)\right)^2,
\end{align*}
%\begin{align*}
%\tau_{\rm R}((FW,g_F);\rho)&=\prod_{q=0}^{p-1}  \left(\frac{\G_q}{l_2^{m-2q}}\right)^{(-1)^q\frac{ r_q}{2}}D_q^{(-1)^q} \prod_{q=p}^{m}  \left(\frac{\G_q}{l_2^{m-2q}}\right)^{(-1)^q\frac{ r_q}{2}}D_q^{(-1)^q} \\
%&=\tau (FW_{(p-1)};\rho,\hat\x) \prod_{q=p}^{m}  \left(\frac{\G_q}{l_2^{m-2q}}\right)^{(-1)^q\frac{ r_q}{2}}
%\hspace{-10pt}[\det\rho( \b_{q+1}( b_{q+1})\semicolon  \zeta_q\semicolon  b_q/ c_q)]\\
%&=\tau (FW_{(p-1)};\rho,\hat\x) \left(\prod_{q=p}^{m}  \left(\frac{\G_q}{l_2^{m-2q}}\right)^{(-1)^q\frac{ r_q}{2}}\right)
%\tau (W_{(p-1)};\rho,\x),
%\end{align*}
where the $\zeta_q$ projects onto an homology basis in $\B$, whose volume element is $\x$. 
This suggests to consider the factor
\begin{align*}
\log&\frac{\tau_{\rm R}((\fr W,g_F);\rho)}{\tau (W_{(p-1)};\rho,\x)\prod_{q=p}^{m}  \left(\frac{\G_q}{l_2^{m-2q}}\right)^{(-1)^q\frac{ r_q}{2}} }=\log\frac{\tau_{\rm R}((W,l_2^2 g);\rho)}{\tau (W_{(p-1)};\rho,\x) }\\
& \hspace{40pt}-\frac{1}{2}\sum_{q=p}^{m} \log l_2^{(-1)^q r_q (m-2q)}+
\frac{1}{2}\sum_{q=0}^{p-1} (-1)^q r_q \log\frac{l_2^{m+1-2q}-l_1^{m+1-2q}}{(m+1-2q)l_2^{m+1-2q}}.
\end{align*}

Hence
\begin{align*}
\log&\frac{\tau_{\rm R}((\fr W,g_F);\rho)}{\tau (W_{(p-1)};\rho,\x)\prod_{q=p}^{m}  \left(\frac{\G_q}{l_2^{m-2q}}\right)^{(-1)^q\frac{ r_q}{2}}
 }=\log\frac{\tau_{\rm R}((W, g);\rho)}{\tau (W_{(p-1)};\rho,\x) }\\
&\hspace{40pt} +\frac{1}{2}\sum_{q=0}^{p-1} \log l_2^{(-1)^q r_q (m-2q)}+
\frac{1}{2}\sum_{q=0}^{p-1} (-1)^q r_q \log\frac{l_2^{2p-2q}-l_1^{2p-2q}}{(2p-2q)l_2^{2p-2q}}.
\end{align*}

For simplicity, we call the above fraction the {\it geometrically regularized} R {\it torsion} of $\fr W$, and we use the notation
\[
\Upsilon_{\rm R}((\fr W,g_F);\rho)=\frac{\tau_{\rm R}((\fr W,g_F);\rho)}{\tau (W_{(p-1)};\rho,\x)\prod_{q=p}^{m}  \left(\frac{\G_q}{l_2^{m-2q}}\right)^{(-1)^q\frac{ r_q}{2}} }.
\]

It is easy to see that the limit is 
\begin{align*}
\lim_{l_1\to 0^+}\log \Upsilon_{\rm R}((\fr W,g_F);\rho_0)=&\frac{1}{2}\log\tau_{\rm R}((W, g_0);\rho)+\frac{1}{2}\sum_{q=0}^{p-1} (-1)^q r_q (2p-1-2q)\log l_2\\
&+\frac{1}{2}\sum_{q=0}^{p-1} (-1)^{q+1} r_q \log\frac{2(p-q)}{l_2},
%\\
%=&\frac{1}{2}\log \tau_{\rm R}((W,l_2^2 g);\rho)+\frac{1}{2}\sum_{q=0}^{p-1} (-1)^{q}  \log\left(\frac{l_2}{2p-2q}\right)^{r_{q}}.
\end{align*}
%\begin{align*}
%\lim_{l_1\to 0^+}\log&\frac{\tau_{\rm R}((\fr W,g_F);\rho)}{\prod_{q=p}^{m}  \left(\frac{\G_q}{l_2^{m-2q}}\right)^{(-1)^q\frac{ r_q}{2}}\tau (W_{(p-1)};\rho,\x) }=\frac{1}{2}\log\tau_{\rm R}((W, g);\rho)\\
%&\hspace{20pt}+\frac{1}{2}\sum_{q=0}^{p-1} (-1)^q r_q (2p-1-2q)\log l_2+\frac{1}{2}\sum_{q=0}^{p-1} (-1)^{q+1} r_q \log\frac{2(p-q)}{l_2}.
%\end{align*}
and this coincides  precisely with  the analytic torsion of the cone up to the boundary term. Comparing with Proposition 4.1  of \cite{HS2}, we see that 
\begin{align*}
\lim_{l_1\to 0^+}\log&\frac{\tau_{\rm R}((\fr W,g_F);\rho_0)}{\tau (W_{(p-1)};\rho_0,\x)\prod_{q=p}^{m}  \left(\frac{\G_q}{l_2^{m-2q}}\right)^{(-1)^q\frac{ r_q}{2}} }=\log I \tau_{\rm R}(C_{l_2} W),
\end{align*}
where the right end side is the intersection R torsion of the cone \cite{Dar1} \cite{HS3}. We have proved an analytic definition of the intersection torsion of a cone, namely:

\begin{theo} In the limit $F W\to C W$, the geometrically  regularized R torsion of the conical frustum $F W$ over an oriented compact connected odd dimensional manifold $W$ gives the Intersection torsion of the cone $C W$ over $W$.
\end{theo}

In order to complete the analysis of the limit case, we need to consider the boundary term. From Section \ref{s4} the boundary term of the frustum is
\[
A_{\rm BM,abs}(\b \fr W)=\frac{1}{2}\int_{W_2} B-\frac{1}{2}\int_{W_1} B,
\]
and hence there is a jump discontinuity with 
\[
\lim_{l_1\to 0^+}A_{\rm BM,abs}(\b \fr W)=\frac{1}{2}\int_{W_2} B=A_{\rm BM,abs}(\b C_{l_2} W),
\]
this completes the proof of the following result.

\begin{theo} In the limit $F W\to C W$, the geometrically  regularized analytic torsion of the conical frustum $F W$ over an oriented compact connected odd dimensional manifold $W$ gives the analytic torsion of the cone $C W$ over $W$.
\end{theo}

\section{The case of a circle} 
\label{lastsec}

Set $\fr =\fr S^1_{\sin\alpha}$, where $S^1_{\sin\alpha}$ denotes a circle of radius $\sin\alpha$. $\fr$  is the finite surface in $\R^3$ parametrized by 
\begin{equation*}\label{tronco do cone} \an =\left\{
\begin{array}{rcl}
x_1&=&r\sin\alpha  \cos{\theta} \\[8pt]
x_2&=&r\sin\alpha  \sin{\theta} \\[8pt]
x_3&=&r\cos\alpha                 \\[8pt]
\end{array}
\right.\end{equation*} 
with $(r,\theta) \in [l_1,l_2]\times[0,2\pi]$, and the metric is the metric induced by the immersion $
g_\an = dr \otimes dr + (\sin^{2}\alpha) r^2 d\theta \otimes d\theta$.

\begin{figure}[h]\centering
 \includegraphics[scale=0.4]{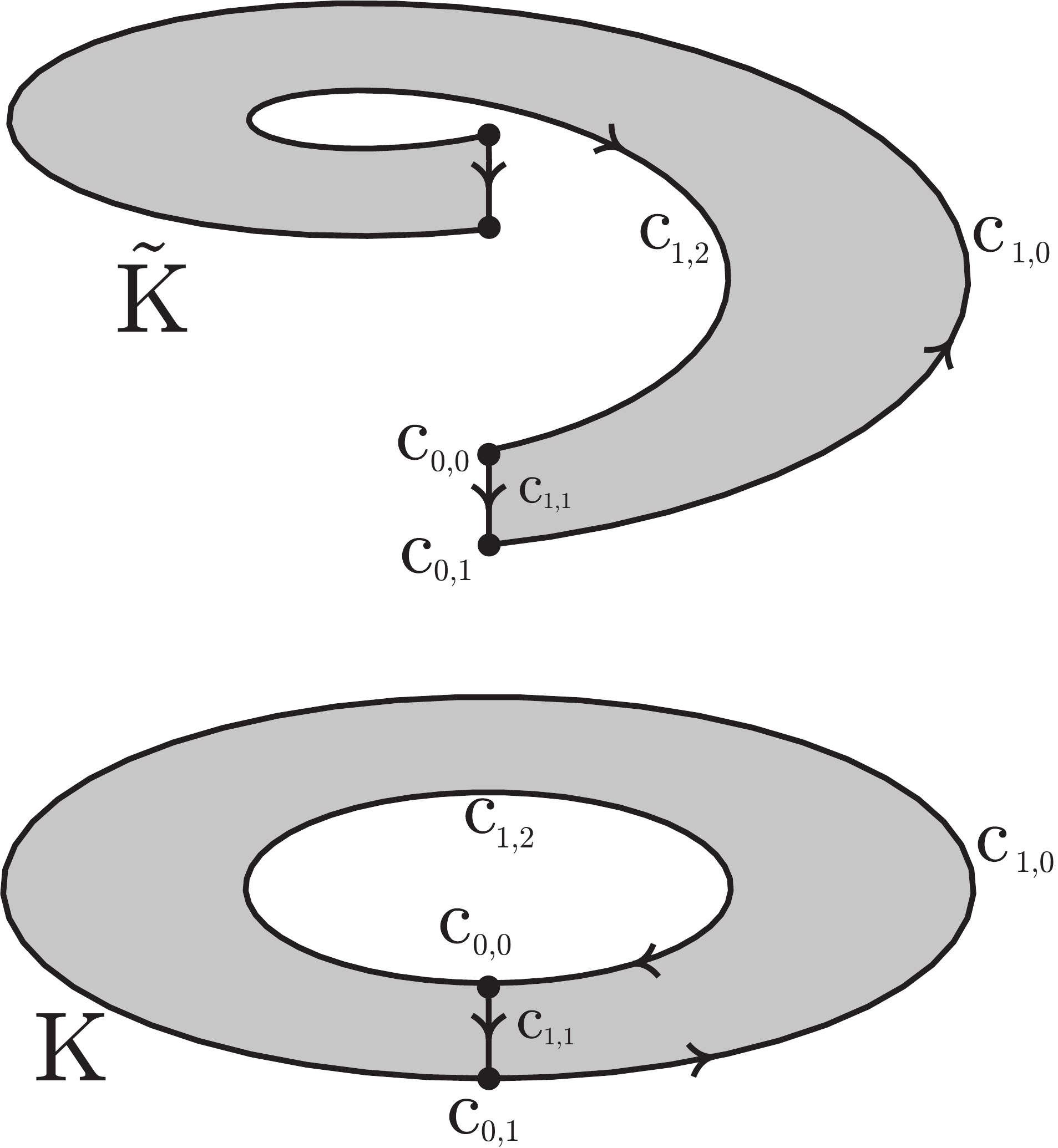}
  \caption{Cell decomposition of $\an$.}\label{fig1}
\end{figure}

Let $K$ denotes the cellular decomposition of $\fr$ described in Figure \ref{fig1}, with subcomplexes $L_j$ that are cellular decompositions of $W_j$. Let $\tilde K$ be the universal covering complex (that is a cellular decomposition of the space $\tilde \fr=[l_1,l_2]\times \R$). It is easy to see that the integral cellular chain complex of $\tilde K$ with the fundamental group acting by covering transformations gives the following  chain complex of $\Z\pi$-modules, where $\pi=\pi_1(\fr)=\Z$,

\[
\xymatrix@C=0.65cm{
\CC(\fr;\Z\pi):&0\ar[r]&\Z\pi[c_2]\ar[r]&\Z\pi[c_{1,0},c_{1,1},c_{1,2}]\ar[r]&\Z\pi[c_{0,0},c_{0,1}]\ar[r]&0,
}
\]
with boundaries
\begin{align*}
\partial_2 (c_{2}) &= c_{0,1} + c_{1,0} - t c_{1,1} - c_{1,2}, \\
\partial_1 (c_{1,0}) &= t c_{0,1} - c_{0,1}, \hspace{1cm}
\partial_1 (c_{1,1}) &=   c_{0,1} - c_{0,0}, \hspace{1cm}
\partial_1 (c_{1,2}) &= t c_{0,0} - c_{0,0}. 
\end{align*}

Taking the trivial representation $\rho_0:\pi\to O(1,\R)$, and  considering the complex of vector spaces $\CC(\fr;\R_\rho)=\R\times_\rho \CC(\fr;\Z\pi)$, we have $H_0(\fr;\R_\rho)=H_1(\fr;\R_\rho)=\R$, $H_2(\fr;\R_\rho)=0$. In order to compute the R torsion, we fix bases for homology. In dimension zero, take  $1 \in \H^{0}(F)$. Then, $\| 1 \| = \sqrt{{\rm Vol} (\an)} = \sqrt{(l^2_2 -l^2_1)\pi\sin\alpha}$, and using the De Rham map $\A_0^{\rm abs} = \P^{-1}_{0} \A^{2}_{\rm rel} \star : \H^{0}(\an)
\rightarrow C_0(K,E_{\rho_0})$ we get 
\begin{align*}
A_0^{\rm abs} \left(\frac{1}{\|1\|}\right) &= \P^{-1}_{0} \A^{2}_{\rm rel} \left(\frac{r \sin\alpha dr \wedge d\theta}{\| 1 \|} \right) = \P^{-1}_{0} (\sqrt{{\rm Vol}(\an)} (c^\dagger_{0,1} + c^\dagger_{0,0})^{*}) \\
&= \sqrt{{\rm Vol}(\an)} (c_{0,1} + c_{0,0}).
\end{align*}

In dimension one,  consider $d\theta \in \H^{1}(\an)$, satisfying absolute BC. We have 
$\| d\theta \|^2 = \frac{2\pi}{\sin\alpha} \ln\frac{l_2}{l_1}$, and we want to apply the De Rham map $\A_{1}^{\rm abs} = - \P^{-1}_{1} \A^{1}_{\rm rel} \star$. Now,  $\A^{1}_{\rm rel}: \H^{1}(\an) \rightarrow C^{1}( K^\dagger, E_\rho)$ is defined by 
$\A^{1}_{\rm rel}(h)( c^\dagger) = \int_{ c^\dagger} h$. Using the basis of $C^{1}( K^\dagger, E_\rho)$, we have $
\A^{1}_{\rm rel}\left( dr\right)(c^\dagger_{1,1}) = 0$, since  $(c^\dagger)^{1}_{1}$ is a circle with constant $r$, and
\begin{align*}
\A^{1}_{\rm rel}\left(-\frac{1}{r\sin\alpha\|d\theta\|} dr\right)( c^\dagger_{1,0} -  c^\dagger_{1,2}) 
&= \int^{l_2}_{l_1} -\frac{1}{r\sin\alpha\|d\theta\|} dr=- \frac{(\ln \frac{l_2}{l_1})^{\frac{1}{2}}}{(2\pi\sin\alpha)^{\frac{1}{2}}}.
\end{align*}

This gives $\A_{1}^{\rm abs}\left(\frac{d\theta}{\|d\theta\|}\right) = \frac{(\ln \frac{l_2}{l_1})^{\frac{1}{2}}}{(2\pi\sin\alpha)^{\frac{1}{2}}} (c_{1,0} -
c_{1,2})$. Therefore: %Fixing the preferred bases $c_q$ for the chain modules and taking the following  bases $b_q$ for boundary and $h_q$ for homology
\begin{align*}
\{\partial_1(b_1),z_0,b_0\} &= \{c_{0,1} - c_{0,0},\emptyset,\sqrt{{\rm Vol}(\an)}(c_{0,0} + c_{0,1})\}, \\
\{\partial_2(b_2),z_1,b_1,\} &= \{c_{1,0} + c_{1,2},c_{1,1},\frac{(\ln
\frac{l_2}{l_1})^{\frac{1}{2}}}{(2\pi\sin\alpha)^{\frac{1}{2}}}(c_{1,0} - c_{1,2})\}, \\
\{\partial_3(b_3),z_2, b_2\} &= \{\emptyset,c_{2},\emptyset\},
\end{align*}
and the torsion is 
\[
\tau_{\rm R}((\an,g_\an);\rho_0) = \prod_{q=1}^2 |\det (\b_{q+1}( b_{q+1}), z_q,b_q/c_q)|=
\frac{\sqrt{{\rm Vol}(\an)}}{\frac{(\ln
\frac{l_2}{l_1})^{\frac{1}{2}}}{(2\pi\sin\alpha)^{\frac{1}{2}}}} = \frac{\pi\sin\alpha(2(l_2^2 - l_1^2))^{\frac{1}{2}}}{(\ln
\frac{l_2}{l_1})^{\frac{1}{2}}}.
\]

Similar calculations for the pairs $(\fr,\b\fr)$ and $(\fr, W_2)$ give, respectively,
\begin{align*}
\tau_{\rm R}((\an,\b \fr ,g_\an);\rho_0) &= \frac{\frac{1}{\sqrt{{\rm
Vol}(\an)}}}{\frac{(2\pi\sin\alpha)^{\frac{1}{2}}}{(\ln\frac{l_2}{l_1})^{\frac{1}{2}}}}= 
\frac{(\ln\frac{l_2}{l_1})^{\frac{1}{2}}}{\pi\sin\alpha (2(l_{2}^{2}-l^{2}_{1}))^{\frac{1}{2}}};\\
\tau_{\rm R}((\an,W_2 ,g_\an);\rho_0) &=1.\\
\end{align*}

Triviality of the last is expected since this corresponds to the cone relative to a point, up to simple homotopy type.

\subsection{Anomaly boundary term} We determine the forms $\mathcal{S}$ and $B$ appearing in the definition of $A_{\rm BS, abs}(\b \fr)$. Since the last is local, and the boundary $\b\fr=W_1\sqcup W_2$ is non connected, we consider the three metrics: $g_\fr$, and $g_j=dr\otimes dr+\sin^2\alpha l_j^2 d\theta\otimes d \theta$. In the first metric, an orthonormal basis is ${e}_r = \frac{\partial}{\partial r}$, ${e}^r = dr$,
    $ {e}_{\theta_1} = \frac{1}{r\sin\alpha}\frac{\partial}{\partial{\theta}}$, and $ {e}^{\theta} = r\sin\alpha d\theta$. 
The non trivial Christoffel symbols are: $\tensor{\Gamma}{ _{\theta} ^r _{\theta}} = -\frac{1}{r}$ and 
$\tensor{\Gamma}{ _{\theta}^{\theta} _r} = \frac{1}{r}$, and the connection one form is
\[
\tensor{(\omega)}{^{r}_{\theta}} = -\frac{1}{r} e^{\theta} = -\sin\alpha d\theta,
\]
implying the vanishing of the curvature. In the metric $g_j$ it is easy to see that the connection one form vanishes identically. Applying the definition 
\begin{align*}
-\mathcal{S}_{W_1} =\mathcal{S}_{W_2}&=  \frac{1}{2} \tensor{(i_j^*\omega_1-i_j^*\omega_{0,l_1})}{^{r}_{\theta}}\hat e^{\theta} =
- \frac{1}{2} \sin\alpha d\theta \hat e^{\theta},
\end{align*}
(where $i_j^*$ is the inclusion of the boundary $W_j$), giving
\begin{align*}
B_{W_j}(\nabla^{T\fr}) =  \frac{1}{2} \int^{1}_{0} \int^{B} \frac{1}{\Gamma(1+\frac{1}{2})} \mathcal{S}_{W_j} du =  \frac{1}{4} \frac{1}{\Gamma(1+\frac{1}{2})} \int^{B} \sin\alpha d\theta \hat e^{\theta}
=\frac{\sin\alpha}{2\pi} d\theta,
\end{align*}
and hence
\begin{align*}
A_{\rm BM,abs}(W_j)%=\log \frac{T_{\rm abs}((\fr,g_\fr);\rho)}{T_{\rm abs}((\fr,g_0);\rho)} 
&= \frac{1}{2}{\rm rk}(\rho) \int_{S^{1}_{l_j\sin\alpha}} B_{W_j}(\nabla_1^{T\fr})  = \frac{(-1)^j}{2}{\rm rk}(\rho)\sin\alpha,\\
A_{\rm BM,rel}(W_j)&= \frac{(-1)^{j+1}}{2}{\rm rk}(\rho)\sin\alpha.
\end{align*}

\subsection{Analytic torsion} In this section we compute the analytic torsion of $\fr$ in the trivial representation. The technique is the one used in  \cite{HS1}, and we refer to that work for details. We write $\frac{1}{\nu}=\sin\alpha$, for convenience. So
%$g_\fr = dr \otimes dr + \frac{r^2}{\nu^{2}} d\theta \otimes d\theta$, 
the Hodge operator is $\star1=\frac{r}{\nu}  d r\wedge d\theta$, $\star d r=  \frac{r}{\nu}   d\theta$, $\star d\theta = -\frac{\nu}{ r}  d r$, 
$\star d r\wedge d\theta= \frac{\nu}{ r}$, and the Laplace operator reads
\begin{align*}
\Delta^{(0)}(f)=&-\b_r^2 f-\frac{1}{r}\b_r
f-\frac{\nu^2}{r^2}\b_{\theta}^2 f,\\
\Delta^{(1)}(f_r d r+ f_{\theta} d_{\theta})=& \left(-\b_r^2 f_r -\frac{\nu^2}{ r^2} \b^2_{\theta} f_r
+\frac{1}{r^2}
f_r-\frac{1}{r}\b_r f_r + \frac{2\nu^2}{r^3}\b_{\theta} f_{\theta} \right)d r\\
&+\left(-\b^2_r f_{\theta} -\frac{\nu^2}{r^2}\b_{\theta}^2 f_{\theta} +\frac{1}{r}\b_r f_{\theta}
-\frac{2}{r}\b_{\theta} f_r\right)
d\theta,\\
\Delta^{(2)}(f dr\wedge d\theta)=&-\b_r^2 f+\frac{1}{r}\b_r f -\frac{2\nu^2}{r^2} f -\frac{1}{r^2}\b_{\theta}^2 f.
\end{align*}

Proceeding as in \cite[Section 3]{HS1}  or \cite[Lemma 3]{HMS}, we find a complete system of eigenforms for $\Delta$, and imposing absolute and relative BC respectively, we obtain the  spectrum 
\begin{align*}
\Sp \Delta_{\rm abs}^{(0)}=& \left\{2:\tilde a_{\nu n,k}^2\right\}_{n,k=1}^{\infty}\cup \left\{\tilde a_{0,k}^2\right\}_{k=1}^{\infty}, \\
\Sp \Delta_{\rm abs}^{(1)}=& \left\{2:\tilde a_{\nu n,k}^2\right\}_{n,k=1}^{\infty} \cup \left\{2:a_{\nu n,k}^2\right\}_{n,k=1}^\infty
\cup \left\{\tilde a_{0,k}^2\right\}_{k=1}^{\infty}\cup \left\{a^2_{0,k}\right\}_{k=1}^{\infty} , \\
\Sp \Delta_{\rm abs}^{(2)}=& \left\{2:a_{\nu n,k}^2\right\}_{n,k=1}^\infty\cup \left\{a^2_{0,k}\right\}_{k=1}^{\infty}; \\
\Sp \Delta_{\rm rel}^{(0)}=& \left\{2:a^2_{\nu n,k}\right\}_{n,k=1}^{\infty}\cup \left\{a^2_{0,k}\right\}_{k=1}^{\infty}, \\
\Sp \Delta_{\rm rel}^{(1)}=& \left\{2:a^2_{\nu n,k}\right\}_{n,k=1}^{\infty} \cup \left\{2:\tilde a_{\nu n,k}^2\right\}_{n,k=1}^\infty
\cup \left\{\tilde a_{0,k}^2\right\}_{k=1}^{\infty}\cup \left\{a^2_{0,k}\right\}_{k=1}^{\infty} , \\
\Sp \Delta_{\rm rel}^{(2)}=& \left\{2:\tilde a_{\nu n,k}^2\right\}_{n,k=1}^\infty\cup \left\{\tilde a_{0,k}^2\right\}_{k=1}^{\infty},
\end{align*}
where the $a_{\nu n,k}, \tilde a_{\nu n,k}$ are the zeros of the function $F_{\nu n}(z), \tilde F_{\nu n,k}$, respectively (here $J_{-0}$ is replaced by $Y_0$):
\begin{align*}
F_{\nu n}(z) &= J_{\nu n}(l_2 z)J_{-\nu n}(l_1 z) - J_{\nu n}(l_1 z)J_{-\nu n}(l_2 z),\\
\tilde F_{\nu n}(z) &= J'_{\nu n}(l_1 z)J'_{-\nu n}(l_2 z) - J'_{\nu n}(l_2 z)J'_{-\nu n}(l_1 z).
\end{align*}

The torsion zeta function is
\[
t_{\rm abs/rel}(s)=\frac{1}{2}\sum_{q=1}^2 (-1)^q \zeta(s,\Delta^{(q)}_{\rm abs/rel}),
\]
and using the above description of the spectra, after some simplification, we obtain

\begin{align}
 t_{\rm abs}(s) =-t_{\rm rel}(s)
&= \sum^{\infty}_{n,k=1} a^{-2s}_{\nu n,k} - \sum^{\infty}_{n,k=1} \tilde a_{\nu n,k}^{-2s} + \frac{1}{2}\sum^{\infty}_{k=1}a_{0,k}^{-2s}-\frac{1}{2} \sum^{\infty}_{k=1}\tilde a_{0,k}^{-2s},\label{tabs}.
%\\
 %t_{\rm rel}(s) &=  \sum^{\infty}_{n,k=1} \tilde a_{\nu n,k}^{-2s} - \sum^{\infty}_{n,k=1} a^{-2s}_{\nu n,k} +\frac{1}{2}
%\sum^{\infty}_{k=1}\tilde a_{0,k}^{-2s} -\frac{1}{2}\sum^{\infty}_{k=1}a_{0,k}^{-2s}\label{trel}.
\end{align}

To compute the derivative at zero of the last two functions, $Z(s, S_0)=\sum_{k=1}^\infty a_{0,k}^{-2s}$, $Z(s, \tilde S_0)=\sum_{k=1}^\infty \tilde a_{0,k}^{-2s}$, we use Proposition 2.4 of \cite{Spr4}. For we need the asymptotic expansion for large $\lambda$ of the Gamma functions associated to the sequences $S_0=\{a_{0,k}\}_{k=1}^\infty $ and $\tilde S_0=\{\tilde a_{0,k}\}_{k=1}^\infty $ (see \cite[Sec. 2.1]{Spr4}). Proceeding as in \cite[Section 5.2]{HS1}, we have the product representations (for $-\pi<\arg(z)<\frac{\pi}{2}$) 
\begin{align*}
G_0(z)&=F_0(i z)=\frac{2}{\pi}\log\frac{l_2}{l_1}{\prod_{k=1}^{+\infty}}\left(1+\frac{z^2}{a^2_{0,k}}\right),\\
 \tilde G_0(z)&=\tilde F_0(i z)=\frac{1}{\pi} \left(\frac{l_2^2 - l_1^2}{l_1 l_2}\right){\prod_{k=1}^{+\infty}}\left(1+\frac{z^2}{(\tilde a_{0,k})^2}\right),
\end{align*}
that give
\begin{align*}
\log \Gamma(-\lambda,S_{0})&=-\log\prod_{k=1}^\infty \left(1+\frac{(-\lambda)}{a_{0,k}^2}\right)
=-\log G_{0}(\sqrt{-\lambda}) + \log\frac{2}{ \pi} +\log\log\frac{l_2}{l_1},\\
\log \Gamma(-\lambda,S_{0})&=-\log\prod_{k=1}^\infty \left(1+\frac{(-\lambda)}{(\tilde a_{0,k})^2}\right)=-\log \tilde G_{0}(\sqrt{-\lambda}) - \log \pi +\log\frac{l_2^2 - l_1^2}{l_1 l_2}.
\end{align*}

Using classical expansions for the Bessel functions, we have the expansions for large $\lambda$
\begin{align*}
\log\Gamma(-\lambda, S_0) &\sim \log\sqrt{-\lambda}  +\frac{1}{2} \log l_1 l_2 + \log \log \frac{l_2}{l_1} +\log 2 -  (l_2 - l_1) \sqrt{-\lambda} %+ O\left(\frac{1}{\sqrt{-\lambda}}\right)
,\\
\log\Gamma(-\lambda, \tilde S_0) &\sim\log\sqrt{-\lambda}  - \frac{1}{2} \log l_1 l_2 + \log (l_2^2 - l_1^2)  -  (l_2 - l_1) \sqrt{-\lambda} %+ O\left(\frac{1}{\sqrt{-\lambda}}\right)
,
\end{align*} 
that give
\begin{align*}
Z'(0,S_0) &= -\frac{1}{2} \log l_1 l_2 - \log \log \frac{l_2}{l_1} -\log 2,\\
Z'_0(0,\tilde S_0)& = \frac{1}{2} \log l_1 l_2 - \log (l_2^2 - l_1^2).
\end{align*}

For the double series we use Theorem 3 of \cite{HS1} and its corollary. For we consider the zeta functions:
$Z(s,S) = \sum^{\infty}_{n,k=1} a^{-2s}_{\nu n,k}$, and $ Z(s,\tilde S) = \sum^{\infty}_{n,k=1} \tilde a_{\nu n,k}^{-2s}$, associated to the double series $S=\{a^{-2s}_{\nu n,k}\}_{n,k=1}^\infty$, $\tilde S=\{\tilde a^{-2s}_{\nu n,k}\}_{n,k=1}^\infty$, respectively. We first prove that the two double sequences are spectrally decomposable over the sequence $U=\{(\nu n)^2\}_{n=1}^\infty$, with 
power 2 and length 3 according to Definition 1 of \cite{HS1} (for the proof see \cite[Section 5.5]{HMS}). Next, we need uniform asymptotic expansion of the associated Gamma functions $\Gamma(\lambda, S_n/(\nu n)^2)$, $\Gamma(\lambda, \tilde S_n/(\nu n)^2)$, for large $n$, and expansions for large $\lambda$. We have  the product representations ($-\pi<\arg(z)<\frac{\pi}{2}$) \cite[Section 5.2]{HS1}
\begin{align*}G_\nu(z)&=F_\nu(i z)=\frac{\sin\pi\nu}{\nu\pi}\left(\frac{l_2^{\nu}}{l_1^{\nu}} -
\frac{l_1^{\nu}}{l_2^{\nu}}\right){\prod_{k=1}^{+\infty}}\left(1+\frac{z^2}{a^2_{\nu,k}}\right),\\
  \tilde G_\nu(z)&=   i^2 \tilde F_\nu(i z)=\frac{\nu\sin\pi\nu}{\pi}\left(\frac{l_2^{\nu}}{l_1^{\nu}} - \frac{l_1^{\nu}}{l_2^{\nu}}\right)\frac{1}{l_1 l_2
z^2}{\prod_{k=1}^{+\infty}}\left(1+\frac{z^2}{\tilde a_{\nu,k}^2}\right).
\end{align*}

This gives
\begin{align*}
\log \Gamma(-\lambda,S_{n}/(\nu n)^2) =& -\log G_{\nu n}(\nu n \sqrt{-\lambda}) + \log\frac{\sin(\nu n\pi)}{\nu n \pi} +\log\left(\frac{l_2^{\nu n}}{l_1^{\nu n}}-\frac{l_1^{\nu n}}{l_2^{\nu n}}\right),\\
\log \Gamma(-\lambda,\tilde S_{n}/(\nu n)^2)=&-\log \tilde G_{\nu n}(\nu n \sqrt{-\lambda}) + \log\frac{\sin(\nu n\pi)}{\nu n\pi}+\log\left(\frac{l_2^{\nu n}}{l_1^{\nu n}}-\frac{l_1^{\nu n}}{l_2^{\nu n}}\right)\\
&-\log(-\lambda  l_1 l_2),
\end{align*}
and using the asymptotic  expansions of  the Bessel functions and of their derivatives for large index \cite[(7.18), Ex. 7.2]{Olv} we have the desired expansions. According to equation (\ref{tabs}), we just need to work with the difference $Z(s,S)-Z(s,\tilde S)$. After some computations, we obtain the asymptotic expansion for large $n$ (uniformly in $\lambda$)
\begin{align*}
\log \Gamma(-\lambda,\tilde S_{n}/(\nu n)^2)  -& \log \Gamma(-\lambda,S_{n}/(\nu n)^2) \\
&= - \frac{1}{2} \log(1-\lambda l_1^2 )(1 -\lambda l_2^2 )+ \phi_1(\lambda)\frac{1}{\nu n} + O\left(\frac{1}{(\nu n)^{2}}\right),
\end{align*}  
where  $\phi_1(\lambda) = V_1(l_1 \sqrt{-\lambda}) -U_1(l_1\sqrt{-\lambda}) - V_1(l_2\sqrt{-\lambda}) + U_1(l_2 \sqrt{-\lambda})$, and using the coefficients in the expansions of the Bessel functions given in \cite[(7.18), Ex. 7.2]{Olv} 
\[
\phi_1(\lambda) = \frac{1}{2}\left(\frac{1}{(1-l_2^2\lambda)^{\frac{1}{2}}}-\frac{1}{(1-l_1^2\lambda)^{\frac{1}{2}}}\right) -\frac{1}{2}\left( \frac{1}{(1-l_2^2\lambda)^{\frac{3}{2}}} - \frac{1}{(1-l_1^2\lambda)^{\frac{3}{2}}}\right).
\]

Applying the definition \cite[(11)]{HS1} 
\begin{align*}
\Phi_{1}(s)&=\int_0^\infty t^{s-1}\frac{1}{2\pi i}\int_{\Lambda_{\theta,c}}\frac{\e^{-\lambda t}}{-\lambda} \phi_{1}(\lambda) d\lambda dt\\
&=(l_2^{2s}-l_1^{2s})\left(\frac{1}{2}\frac{\Gamma(s+\frac{1}{2})}{\Gamma(\frac{1}{2})s} - \frac{1}{2}\frac{\Gamma(s+\frac{3}{2})}{\Gamma(\frac{3}{2})s}\right),
\end{align*} 
and $\Rz_{s=0} \Phi_1(s) = \Ru_{s=0}\Phi_1(s) = 0$. On the other side, the expansion for large $\lambda$ is 
\begin{align*}
&\log\Gamma(-\lambda, \tilde S_n/(\nu n)^2) - \log\Gamma(-\lambda, S_n/(\nu n)^2) = -\log (-\lambda) -\log l_1 l_2 + O\left(\frac{1}{\sqrt{-\lambda}}\right).
\end{align*}

Recalling the definition in \cite[(13)]{HS1}, we have $a_{0,0,n}=-\log l_1 l_2$, $a_{0,1,n}=-1$, $b_{1,0,0}=b_{1,0,1}=0$, that gives
\begin{align*}
A_{0,0}(s)&=-\sum_{n=1}^\infty \frac{\log l_1 l_2}{(\nu n)^{2 s}} ,&
A_{0,1}(s)&=-\sum_{n=1}^{\infty}\frac{1}{(\nu n)^{2s}},
\end{align*}  and 
\begin{align*}
A_{0,0}(0) &=\frac{1}{2}\log l_1 l_2, &
A_{0,1}(0) &=\frac{1}{2},&
A'_{0,1}(0) &=-\log\nu + \log2\pi.
\end{align*}

Collecting all these results and applying Theorem 3 of \cite{HS1} we obtain
\[
Z'(0,\tilde S) - Z'(0,S) = -\frac{1}{2}\log l_1 l_2 + \log\nu -\log2\pi,
\]
that by equation (\ref{tabs}) gives
\begin{align*}
\log T_{\rm abs}((\fr,g_\fr);\rho_0)=-\log T_{\rm rel}((\fr,g_\fr);\rho_0)=t'_{\rm abs}(0)  
=\log\frac{\pi (2(l_2^2 - l_1^2))^{\frac{1}{2}}}{(\log \frac{l_2}{l_1})^{\frac{1}{2}}\nu}.
\end{align*}

\subsection{Some limits}\label{last} The geometry of $\fr$ has at least two natural limit cases: the cone over $W_2$, reached for $l_1\to 0^+$, and the cylinder over $W_2$, reached for $l_1=l_2$.  We investigate in this section the value of the torsion of $\fr$ in these two limit cases. The first case is an instance of the general case discussed in Section \ref{lim}. It is easy to realize that in the limit for $l_1\to 0^+$, the  torsion of the conical frustum diverges. So consider the geometric regularized R torsion (the extension to analytic torsion is straightforward). Then, $W=S^1_{l_2\sin\alpha }$, the circle of radius $l_2\sin\alpha $, and $W_{(0)}$ is its preferred 0 cell. Since $m=1=p$, the volume element $\x$ is fixed by the Ray and Singer basis for the homology of $S^1$. Thus,
\[
\tau_{\rm R}((S^1_{l_2\sin\alpha })_{(0)};\rho_0, \x)=[\det\rho(h_0/c_0)].
\]

Looking at the calculation of the torsion of a sphere  in \cite{MS}, the harmonic basis in dimension 0 is $1/\sqrt{{\rm Vol}(S^1_{l_2\sin\alpha })}$, and applying the De Rham map $h_0=\left\{\sqrt{{\rm Vol}(S^1_{l_2\sin\alpha })}c_{0,0}\right\}$. Since $\rho=\rho_0$, then $\tau_{\rm R}((S^1_{l_2\sin\alpha })_{(0)};\rho, \x)=\sqrt{{\rm Vol}(S^1_{l_2\sin\alpha })}=\sqrt{2\pi l_2\sin\alpha}$. Therefore,
\[
\Upsilon_{\rm R}((\fr,g_F);\rho_0)= \frac{\tau_{\rm R}((\fr,g_F);\rho_0)}{\tau_{\rm R} ((S^1_{l_2 \sin\alpha})_{(0)};\rho_0,\x)\left(\frac{\G_1}{l_2^{-1}}\right)^{-\frac{ r_1}{2}}} =\frac{\tau_{\rm R}((\fr,g_F);\rho_0)}{\sqrt{\frac{\pi(l_2^2-l_1^2)\sin\alpha}{l_2\log\frac{l_2}{l_1}}}}.
\]

This gives
\[
\lim_{l_1\to 0^+} \Upsilon_{\rm R}((\fr,g_F);\rho_0)=\sqrt{{\rm Vol}(C S^1_{l_2\sin\alpha})}=\sqrt{\pi l^2_2\sin\alpha}=\tau_{R}((C S^1_{l_2\sin\alpha},g_C);\rho_0).
\]

Next, consider the cylinder. This is reached by fixing $b_1=l_1\sin\alpha$ and $h=(l_2-l_1)\cos\alpha$,  and taking the limit for $\alpha\to 0$. The result is
\[
\lim_{b_1,h~{\rm fixed}, \alpha\to 0^+} T_{\rm abs}((\fr,g_\fr);\rho_0)=2\pi b_1={\rm Vol}(S^1_{b_1})=\tau_{\rm R}((S^1_{b_1}, g_{S^1_{b_1}});\rho_0),
\]
where the $g_{S^1_{b_1}}$ is the standard metric (see \cite{MS}), consistently with  the fact that a cylinder has the same simple homotopy type of a circle.


\begin{thebibliography}{99}

\bibitem{BM1} J. Br\"uning and Xiaonan Ma, {\em An anomaly formula for Ray-Singer metrics on manifolds with boundary}, GAFA 16 (2006) 767-873.


\bibitem{BM2} J. Br\"uning and Xiaonan Ma, {\em On the gluing formula for the analytic torsion}, to appear.

\bibitem{Che0} J. Cheeger, {\em On the spectral geometry of spaces with conical singularities}, Proc. Nat. Acad. Sci. 76 (1979) 2103-2106.

\bibitem{Che1} J. Cheeger, {\em Analytic torsion and the heat equation}, Ann. Math. 109 (1979) 259-322.

\bibitem{DF} X. Dai and H. Fang {\em Analytic torsion and R-torsion for manifolds with boundary},  Asian J. Math. 4 (2000) 695-714.


\bibitem{Dar1} A. Dar, {\em Intersection R-torsion and the analytic torsion for pseudomanifolds}, Math. Z. 154 (1987) 155-210.

%\bibitem{GZ} I. S. Gradshteyn and I. M. Ryzhik, {\em Table of integrals, Series and Products}, Academic Press, 2007.

\bibitem{HMS} L. Hartmann, T. de Melo and M. Spreafico, {\em Reidemeister torsion and analytic torsion of discs}, Bollettino U.M.I. 9 (2009) 529-533, arXiv:0811.3196v1.



\bibitem{HS1} L. Hartmann and M. Spreafico, {\em The analytic torsion of a cone over a sphere}, J. Math. Pure Ap. 93 (2010) 408-435.


\bibitem{HS2} L. Hartmann and M. Spreafico, {\it The Analytic Torsion of the Cone over an Odd Dimensional Manifold}, J. Geom. Phys. 61 (2011) 624-657.


\bibitem{HS3} L. Hartmann and M. Spreafico, {\it An extension of the Cheeger-M\"{u}ller theorem for a  cone}, 

\bibitem{MS} T. de Melo and M. Spreafico, {\it Reidemeister torsion and analytic torsion of spheres} J.  Homotopy Rel. Struc., 4 (2009)181-185.

\bibitem{Mil} J. Milnor, {\em Whitehead torsion}, Bull. AMS 72 (1966) 358-426.

\bibitem{Mul} W. M\"{u}ller, {\em Analytic torsion and R-torsion of Riemannian manifolds}, Adv. Math. 28 (1978) 233-305.


\bibitem{Mul2} W. M\"{u}ller, {\it Analytic torsion and R-torsion for unimodular representations}, J. Amer. Math. Soc. 6 (1993) 721-753.


\bibitem{Luc} W. L\"uck,  \emph{Analytic and topological torsion for manifolds with boundary and symmetry},  J. Differential Geom.  37   (1993) 263-322.

\bibitem{Nag} M. Nagase, {\em De Rham-Hodge theory on a manifold with cone-like singularities}, Kodai Math. J., 1 (1982) 38-64.

\bibitem{Olv} F.W.J. Olver, {\em Asymptotics and special functions}, AKP, 1997.

\bibitem{RS} D.B. Ray and I.M. Singer, {\em R-torsion and the Laplacian
on Riemannian manifolds}, Adv. Math. 7 (1971) 145-210.


\bibitem{Spr4} M. Spreafico, {\em Zeta invariants for sequences of spectral type, special functions and the Lerch formula}, Proc. Roy. Soc. Edinburgh 136A (2006) 863-887.




\end{thebibliography}
\end{document}